\documentclass[12pt]{article}
\usepackage[mathscr]{eucal}
\usepackage{amssymb,amsmath}
\font\cmc=cmcsc10  scaled \magstep2

\newcommand\lag{\langle}
\newcommand\rag{\rangle}

\newcommand\R{\mathbb{R}}

\newcommand\vk{\vskip}
\newcommand\hk{\hskip}

\newcommand\al{\alpha}

\newcommand\iy{\infty}

\newcommand\la{\lambda}

\newcommand\el{\ell}

\newcommand\varep{\varepsilon}
\newcommand\rg{\rightarrow}
\newcommand\lrg{\longrightarrow}

\newcommand\De{\Delta}

\newcommand\ov{\overset}
\newcommand\und{\underset}
\newcommand\no{\noindent}

\newcommand\ovl{\overline}

\newcommand\col{\colon\hk-.5em}

\newtheorem{Proof.}{\it Proof.}

\begin{document}
\vbox to .5truecm{}

\begin{center}
\cmc A Geometric Characterization of Extremal Sets
\end{center}

\begin{center}
\cmc in {\mathversion{bold} $\el_p$} Spaces
\end{center}
\vk.3cm
\begin{center}
by Viet NguyenKhac
\footnote{Institute of Mathematics, P.O.Box 631 Bo Ho, 
10000 Hanoi, Vietnam;\ e-mail:\ nkviet@thevinh.ncst.ac.vn}  
\ \& Khiem NguyenVan
\footnote{Department of Mathematics and Informatics, Hanoi University of Education, Cau Giay dist., Hanoi, Vietnam.}
\end{center}
\begin{center} 
Hanoi Institute of Mathematics \& Hanoi University of Education
\end{center}
\vk.1cm
\begin{center}
1991 Mathematics Subject Classification. Primary 46B20, 46E30.
\end{center}

\normalsize
\begin{abstract} {We give a geometric characterization of extremal sets in $\el_p$ spaces that generalizes our previous result for such sets in Hilbert spaces.}
\end{abstract}
\vk.5cm					
\begin{center} 
\textbf {1.\ \ Introduction}
\end{center}
\vk.2cm

Let\ $(X,\|\cdot\|)$\ be a Banach space. For a non-empty bounded subset\ $A$\ of\ $X$\ and a non-empty subset\ $B$\ of\ $X$\ let's fix the following notations:\ \ $d(A)\col=$\ sup $\{\|x-y\|:\ \ x, y\in A\}$\ -- the diameter of\ $A$;\ \ $r_B(A)\col=\und{y\in B}{\text {inf}}\ \und{x\in A}{\text {sup}}\ \|x-y\|$\ -- the relative Chebyshev radius of\ $A$\ with respect to\ $B$;\ \ in particular\ $r(A)\col=r_{\ovl{co} A}(A)$ with $\ovl{co} A$\ denoting the closed convex hull of\ $A$;\ a point\ $y\in B$\ is called a Chebyshev center of\ $A$\ in\ $B$,\ \ if\ \ $\und{x\in A}{\text {sup}}\ \|x-y\|=r_B(A)$.  
\vk.2cm
The Jung constant of\ $X$\ is defined by\ $J(X)\col=$ sup $\{r_X(A):\ A\subset X,\ \ {\text {with}}\ \ d(A)=1\}$.\ The problem of estimating Jung's constant plays an important role in the geometry of Banach spaces ({\it cf.} \cite{Am}, \cite{Pich}). It is well-known that for inner-product spaces\ $J(E^n)=\sqrt{\dfrac{n}{2(n+1)}}$\ \ (\cite{Jun}, {\it cf.} \cite{Dgk}, \cite{IoT}) and\ $J(H)=\dfrac{1}{\sqrt{2}}$\ ($H$\ denotes a Hilbert space) (\cite{Rou}). In general if\ $X$\ is an\ $n$-dimensional normed space, then\ $J(X)\le\dfrac{n}{n+1}$.\ Furthermore the equality is attained for certain spaces (see \cite{Pich}). As for\ $\el_p,\ L_p$\ spaces\ $(p>1)$\ S. A. Pichugov (\cite{Pich}) has obtained the exact values for the Jung constant of these spaces:\ \ $J(\el_p)=J(L_p)=\dfrac{1}{\sqrt[q]{2}},\ q\col=\dfrac{p}{p-1}$,\ if $1<p\le 2$, and $J(\el_p)=J(L_p)=\dfrac{1}{\sqrt[p]{2}}$,\ if $p>2$. 
\vk.2cm
\textbf{Definition\ 1.1}\ (\cite{Nkk}).\quad We say that a bounded subset\ $A$\ of\ $X$\ consisting of at least two points is an extremal set, if\ $r_X(A)=J(X).d(A)$.
\vk.2cm
The main result of \cite{Nkk} states that a bounded subset\ $A$\ of a Hilbert space\ $H$\ with\ $r(A)=1$\ is extremal if and only if for every\ $\varep\in (0,\sqrt{2})$,\ for every positive integer\ $m$\ there exists an\ 
$m$-simplex\ $\De$\ with its vertices in\ $A$\ and each edge of\ $\De$\ has length not less than\ $\sqrt{2}-\varep$.\ Furthermore for such a subset\ $A$\ we have\ $\al(A)=\sqrt{2}$\ and $\chi(A)=1$,\ where\ $\al(A),\ \chi(A)$\ denote the Kuratowski and Hausdorff measures of non-compactness of\ $A$\ respectively.
\vk.2cm
Our aim in this paper is to treat the next interesting case:\ \ the case of\ $\el_p$\ spaces\ $(1<p<\iy)$.\ We obtain a partial generalization of the result above. More precisely, if\ $A$\ is an extremal subset of a given\ $\el_p$\ space\  $(1<p<\iy)$,\ then\ $\al(A)=d(A)$.\ As an immediate consequence one obtains a Gulevich-type result for\ $\el_p$ spaces:\ \ extremal sets in\ $\el_p,\ (1<p<\iy)$\ are not relatively compact ({\it cf.}\ \cite{Gul}). Moreover for every\ $\varep\in (0,d(A))$,\ for every positive integer\ $m$\ there exists an\ $m$-simplex\ $\De$\ with its vertices in\ $A$\ and each edge of\ $\De$\ has length not less than\ $d(A)-\varep$.\ The proof is based on a further development of a purely combinatorial method in our previous paper \cite{Nkk} which essentially relies on a very deep part of convex analysis. It should be noted that this observation was first noted in \cite{IoT}, \S 10.2, where the authors exposed  classical Jung's theorem from the point of view of ``subdifferentials", and  was later extended in \cite{Pich} to the case of\ $\el_p$\ spaces.
\vk.2cm
The paper is organized as follows. In \S 2 we review some facts related to the ``Clearance"-type Theorem (``Decomposition Theorem" in translation), especially for\ $\el_p$\ spaces, we shall need in the sequel. The heart of this section is Proposition 2.1 which should be considered as an infinite-dimensional variation of the main claim in \cite{Pich}. In \S 3 we first formulate without proof two auxiliary inequalities (Lemma 3.1) one of which was essentially due to N. I. Chernykh (personal communication to the author of \cite{Pich}). Our main results are Theorems 3.2 and 3.4 the proof of which is proceeded then by means of the combinatorial method of \cite{Nkk} extended to\ $\el_p$\ spaces together with Proposition 2.1 and inequalities from Lemma 3.1.
\vk.5cm
\begin{center}
\textbf {2.\ \ Preliminaries}
\end{center}
\vk.2cm

For the basic definitions and concepts in convex analysis we refer the reader to \cite{IoT}. The following proposition is a slight generalization of the main proposition of \cite{Pich}.
\vk.5cm
\textbf{Proposition\ 2.1.}\quad {\it Let\ $X$\ be a uniformly smooth Banach space,\ $A=\{x_1,\ x_2,\ \ldots,\ x_n\}$\ a finite subset with\ $r\col=r_X(A)>0$,\
$c$\ a Chebyshev center of\ $A$\ in\ $X$.\ Then there exist points\ $y_1,\ y_2,\ \ldots,\ y_m\ (m\le n)$\ in\ $A$,\ linear functionals\ $f_1,\ f_2,\ \ldots,\ f_m$\ in\ $X^*$\ and positive numbers\ $\al_1,\ \al_2,\ \ldots,\ \al_m$\ satisfying:
\vk.2cm
{\rm (i)}\ \ $\lag y_i-c,f_i\rag\col=f_i(y_i-c)=\|y_i-c\|=r$\ for\ \ $i=1,\ 2,\ \ldots m$;
\vk.1cm
{\rm (ii)}\ \ $\|f_i\|_{X^*}=1,\ \ i=1,\ 2,\ \ldots,\ m$;
\vk.1cm
{\rm (iii)}\ \ $\und{i=1}{\ov{m}{\sum}}\ \al_i\ f_i=0,\ \ \und{i=1}{\ov{m}{\sum}}\ \al_i=1$.  }
\vk.2cm
{\it Proof}.\quad For completeness we give here a proof valid in all dimensions. From the uniform smoothness of\ $X$\ it follows that\ $X$\ is reflexive (hence there exists a Chebyshev center of\ $A$\ in\ $X$)\ and the mapping

  $$\begin{aligned} J\colon & X\lrg 2^{X^*}\\
   &\ x\longmapsto J(x)\col=\{ x^*\in X^*:\ \ \|x^*\|=1,\ \lag     x,x^*\rag=\|x\|\}\end{aligned}$$

\no is single-valued. Let's consider the following two functionals:

  $$\begin{aligned} F\colon & X\times A\lrg \R\\
   &\ (x,a)\ \longmapsto F(x,a)\col=\|x-a\|,\end{aligned}$$

\no and

  $$\begin{aligned} f\colon & X\lrg \R\\
   &\ x\longmapsto f(x)\col=\und{a\in A}{\max}\ F(x,a)\end{aligned}.$$

With\ $A_0(x)\col=\{ a\in A:\ \ F(x,a)=f(x)\}$\ we see that the hypotheses of Theorem 3, \S 4.2, Chapter 4 in \cite{IoT} are fulfilled.
Hence for every\ $x\in X$\ we have

  $$\ovl{co}\ \Big( \und{a\in A_0(x)}{\bigcup}\ \partial F(x,a) \Big) =\partial f(x)\eqno{(2.1)}$$   

\no where\ $\partial F(x,\cdot)$\ and\ $\partial f(x)$\ denote subdifferentials of\ $F$\ and\ $f$\ at\ $x$,\ and the closure in (2.1) is taken in the\ $w^*$-topology of the space\ $X^*$,\ which clearly coincides with the\ $w$-topology of\ $X^*$,\ since\ $X$\ is reflexive.
\vk.2cm
Recall that for\ $c\in X$\ to be a Chebyshev center of\ $A$\ in\ $X$\ it is necessary and sufficient that\ $0\in\partial f(c)$\ ({\it cf.} \cite{IoT}, \S 1.3, Proposition 1). One may write\ $A_0(c)=\{y_1,\ y_2,\ \ldots,\ y_m\}$ with $\|y_i-c\|=r$\ \ for\ \ $i=1,\ 2,\ \ldots,\ m$.\ From (2.1) it follows that\  
$0\in \ovl{co}\ \big(\und{i=1}{\ov{m}{\bigcup}}\ \partial F(c,y_i)\big)$.\ Since\ $\|y_i-c\|=r>0$\ we have

  $$\partial F(c,y_i)=\{ x^*\in X^*:\ \ \|x^*\|=1,\ \lag y_i-c,x^*\rag=r\}=J(y_i-c).$$

As noted above\ $J$\ is a single-valued mapping, hence\ $J(y_i-c)$\ consists of a unique point, say\ $f_i$.\ Therefore\ $0\in \ovl{co}\ \{f_1,\ f_2,\ \ldots,\ f_m\}=co\ \{f_1,\ f_2,\ \ldots,\ f_m\}$,\ and so there exist non-negative numbers\ $\al_1,\ \al_2,\ \ldots,\ \al_m$\ such that\ $\und{i=1}{\ov{m}{\sum}}\ \al_i=1$\ and\ $0=\und{i=1}{\ov{m}{\sum}}\ \al_i\ f_i$.\ Without loss of generality one can assume that all\ $\al_i,\ i=1,\ 2,\ \ldots,\ m$\ are positive. It is a simple verification that these data also satisfy the conditions (i), (ii) above. The proof is complete.
\vk.5cm
\textbf{Remark\ 2.2.}\quad We shall be interested mainly in the case\ $X=\el_p\ (p>1)$.\ For this purpose it is more convenient to use the following ``scaled" version of the mapping\ $J$:\ \ for\ $x\in X$\ and\ $p\in (1,\iy)$\ we define\ $J(x)\col=\{ x^*\in X^*:\ \ \lag x,x^*\rag=\|x\|.\|x^*\|=\|x\|^p\}.$\ Obviously for\ $f_i, y_i$\ and\ $c$\ as in Proposition 2.1,\ $\|y_i-c\|^{p-1}.f_i\in J(y_i-c)$.
\vk.5cm
\textbf{Remark\ 2.3.}\quad It is well-known that spaces\ $\el_p\ \ (1<p<\iy)$\ are both uniformly convex and uniformly smooth. Hence the mapping\ $J$\ is single-valued. Also it is weakly sequentially continuous in the following sense:\ if\ $\{x_n\}$\ converges weakly to\ $x$\ in\ $\el_p$,\ then\ $\{J(x_n)\}$\ converges weakly to\ $J(x)$\ in\ $\el_q$\ ($q=\dfrac{p}{p-1}$).
\vk.5cm
\textbf{Remark\ 2.4.}\quad Clearly if\ $x=(x_1,\ x_2,\ \ldots,\ x_n, \ldots )\in\el_p\ \ (1<p<\iy)$,\ then\ $J(x)=($sgn$(x_1).|x_1|^{p-1},$\ sgn$(x_2).|x_2|^{p-1},\ \ldots,\ $sgn$(x_n).|x_n|^{p-1},\ \ldots )$.\ Hence in the situation of Proposition 2.1 with\ $X=\el_p$\ we have

  $$f_i=\dfrac{J(y_i-c)}{{\quad\ }\|y_i-c\|^{p-1}}=\dfrac{J(y_i-c)}{\ \ r^{p-1}}.$$ 
\vk.3cm					
\begin{center} 
\textbf {3.\ \ The results}
\end{center}
\vk.2cm

For our use later it is convenient to formulate two auxiliary inequalities in the following lemma. 
\vk.5cm
\textbf{Lemma\ 3.1.}\quad {\it Let\ $a$\ and\ $b$\ be two real numbers. Then
\vk.2cm
{\rm (i)}\ \ $\Big|${\rm sgn}$(a).|a|^{p-1}-${\rm sgn}$(b).|b|^{p-1}\Big|\le 2^{2-p}.|a-b|^{p-1}$,\ \ provided\ $1<p\le 2$;
\vk.1cm
{\rm (ii)}\ \ $|a|^p+|b|^p-p.\big( a.${\rm sgn}$(b).|b|^{p-1}+b.${\rm sgn}$(a).|a|^{p-1}\big)\le |a-b|^p$,\ \ if $p>2$.}
\vk.4cm
The proof is standard, and we shall omit it.
\vk.5cm
\textbf{Theorem\ 3.2.}\quad {\it Let\ $A$\ be an extremal set in an\ $\el_p$\ space with\ $1<p<\iy$.\ Then we have\ $\al(A)=d(A)$.}
\vk.2cm
Here\ $\al(A)$\ denotes the Kuratowski measure of non-compactness of\ $A$.
\vk.2cm
{\it Proof}.\quad We may assume\ $r_{\el_p}(A)=1$.\ Then for each integer number\ $n\ge 2$\ we have\ $\und{x\in A}{\bigcap}\ B(x,1-\dfrac{1}{n})=\emptyset$,\ where\ $B(x,r)$\ denotes the closed ball centered at\ $x$\ with radius\ $r$\ which is weakly compact since\ $\el_p$\ is reflexive. Hence there exist\  $x_{t_{n-1}+1},\ x_{t_{n-1}+2},\ \ldots,\ x_{t_n}$\ in\ $A$\ such that\ $\und{i=t_{n-1}+1}{\ov{t_n}{\bigcap}}\ B(x_i,1-\dfrac{1}{n})=\emptyset$\ (with convention\ $t_1=0$). 
\vk.2cm
Setting\ $A_n\col=\{x_{t_{n-1}+1},\ x_{t_{n-1}+2},\ \ldots,\ x_{t_n}\}$\ we denote the Chebyshev center of\ $A_n$\ in\ $\el_p$\ by\ $c_n$\ and let\ $r_n\col=r_{\el_p}(A_n)$,\ then\ $r_n>1-\dfrac{1}{n}$.\ In view of Proposition 2.1 one can find\ $y_{s_{n-1}+1},\ y_{s_{n-1}+2},\ \ldots,\ y_{s_n}$\ in\ $A_n$,\ continuous linear functionals\ $f_{s_{n-1}+1},\ f_{s_{n-1}+2},\ \cdots,\ f_{s_n}$\ on\ $\el_p$\ and positive numbers\ $\al_{s_{n-1}+1},\ \al_{s_{n-1}+2},\ \ldots,\ \al_{s_n}$\ (with convention\ $s_1=0$)\ such that:
\vk.2cm
{\rm (i)}\ \ $\lag y_i-c_n,f_i\rag=r_n$\ for\ \ $i=s_{n-1}+1,\ s_{n-1}+2,\ \ldots,\ s_n$;
\vk.1cm
{\rm (ii)}\ \ $\|f_i\|=1,\ \ i=s_{n-1}+1,\ s_{n-1}+2,\ \ldots,\ s_n$;
\vk.1cm
{\rm (iii)}\ \ $\und{i=s_{n-1}+1}{\ov{s_n}{\sum}}\ \al_i\ f_i=0,\ \ \und{i=s_{n-1}+1}{\ov{s_n}{\sum}}\ \al_i=1$.
\vk.2cm
Setting\ $A_\iy\col=\{y_{s_{n-1}+1},\ y_{s_{n-1}+2},\ \ldots,\ y_{s_n}\}_{n=2}^\iy$\ we claim that\ $\al(A_\iy)=d(A)$.\ Suppose on the contrary\ $\al(A_\iy)<d(A)$.\ Then one can choose\ $\varep_0\in (0,d(A))$\ satisfying\ $\al(A_\iy)\le d(A)-\varep_0$,\ and so subsets\ $D_1,\ D_2,\ \cdots,\ D_m$\ of\ $\el_p$\ with\ $d(D_i)\le d(A)-\varep_0$\ for every\ $i=1,\ 2,\ \ldots,\ m$,\ such that\ $A_\iy\subset\und{i=1}{\ov{m}{\bigcup}}\ D_i$.\ There exists at least one set among\ $D_1,\ D_2,\ \cdots,\ D_m$,\ say\ $D_1$\ with the property that there are infinitely many\ $n$\ satisfying

  $$\sum_{i\in J_n}\ \al_i\ge \dfrac{1}{m}\eqno{(3.1)}$$

\no where

  $$J_n\col=\{i\in I_n\col=\{s_{n-1}+1, s_{n-1}+1, \ldots, s_n\}:\ \ y_i\in D_1\}.$$
\vk.2cm
We shall estimate the sum\ $T_n\col=\und{i, j\in I_n}{\sum}\ \al_i.\al_j.\lag y_i-y_j,f_i-f_j\rag$.\ We have

  $$\begin{aligned} T_n & =\und{i, j\in I_n}{\sum}\ \al_i.\al_j.\lag (y_i-   c_n)-(y_j-c_n),f_i-f_j\rag=\\
  & =\und{i, j\in I_n}{\sum}\ \al_i.\al_j.\Big[ \lag y_i-c_n,f_i\rag+
  \lag y_j-c_n,f_j\rag-\lag y_i-c_n,f_j\rag-\lag y_j-c_n,f_i\rag\Big]=\\
  & =2.\und{i, j\in I_n}{\sum}\ \al_i.\al_j.r_n-\und{i\in I_n}{\sum}\       \al_i\big\lag y_i-c_n,\und{j\in I_n}{\sum}\ \al_j.f_j\big\rag-
  \und{j\in I_n}{\sum}\ \al_j\big\lag y_j-c_n,\und{i\in I_n}{\sum}\   \al_i.f_i\big\rag=\\
  & =2.r_n-2.\und{i\in I_n}{\sum}\ \al_i\lag y_i-c_n,0\rag=2.r_n\hk8.7cm {(3.2)}
  \end{aligned}$$ 

As noted in Remark 2.4
   
  $$f_i=\dfrac{J(y_i-c_n)}{\quad\ \|y_i-c_n\|^{p-1}}=\dfrac{J(y_i-c_n)}{\ \ r_n^{p-1}},\ \forall i\in I_n.$$

\no Therefore for all\ $i, j\in I_n$

  $$\lag y_i-y_j,f_i-f_j\rag=\dfrac{1}{r_n^{p-1}}.\big\lag (y_i-c_n)-(y_j-c_n),J(y_i-c_n)-J(y_j-c_n)\big\rag.$$
\vk.2cm
1)\ \ {\it The case}\ $1<p\le 2$:\ \ Since the Jung constant\ $J(\el_p)=\dfrac{1}{\sqrt[q]{2}}$,\ thus\ $d(A)=\sqrt[q]{2}$.\ Applying part (i) of Lemma 3.1 coordinatewise to the above expression and remembering\ $q=\dfrac{p}{p-1}$\ one gets

  $$\big\lag (y_i-c_n)-(y_j-c_n),J(y_i-c_n)-J(y_j-c_n)\big\rag\le 2^{2-p}.\|(y_i-c_n)-(y_j-c_n)\|^p.$$
\vk.2cm
\no Hence

  $$\begin{aligned} T_n & \le\dfrac{1}{r_n^{p-1}} \und{i, j\in I_n}{\sum}\ \al_i.\al_j.2^{2-p}.\|y_i-y_j\|^p=\\
  & =\dfrac{2^{2-p}}{r_n^{p-1}}.\Big(\und{i, j\in J_n}{\sum}\ \al_i.\al_j.\|y_i-y_j\|^p+\und{(i,j)\in I^2_n\setminus J^2_n}{\sum}\ \al_i.\al_j.\|y_i-y_j\|^p\Big)\le\\
  &\le \dfrac{2^{2-p}}{r_n^{p-1}}.\bigg[ \big(\sqrt[q]{2}-\varep_0\big)^p .\und{i, j\in J_n}{\sum}\ \al_i.\al_j+2^{{}^\frac{p}{q}}.\Big(1-\und{i, j\in J_n}{\sum}\ \al_i.\al_j\Big)\bigg]=\\
  & =\dfrac{2^{2-p}.2^{{}^\frac{p}{q}}}{r_n^{p-1}}-\dfrac{2^{2-p}}{r_n^{p-1}}.\Big( 2^{{}^\frac{p}{q}}-\big(\sqrt[q]{2}-\varep_0\big)^p\Big).\Big( 
\und{i, j\in J_n}{\sum}\ \al_i.\al_j\Big)\le\\
  & ({\text {for all $n$ satisfying}}\ (3.1))\\
  &\le \dfrac{2^{2-p}.2^{{}^\frac{p}{q}}}{r_n^{p-1}}-\dfrac{2^{2-p}}{r_n^{p-1}}.\Big[ 2^{{}^\frac{p}{q}}-\big(\sqrt[q]{2}-\varep_0\big)^p\Big].\dfrac{1}{m^2}=\\
  &=\dfrac{2}{r_n^{p-1}}-\dfrac{2^{2-p}}{m^2.r_n^{p-1}}.\Big[ 2^{{}^\frac{p}{q}}-\big(\sqrt[q]{2}-\varep_0\big)^p\Big]\hk8cm {(3.3)}
  \end{aligned}$$
\vk.2cm
Comparing (3.2) and (3.3) we obtain

  $$2.r_n^p\le 2-\dfrac{2^{2-p}}{m^2}.\Big[ 2^{{}^\frac{p}{q}}-\big(\sqrt[q]{2}-\varep_0\big)^p\Big]$$
  
\no for all\ $n$\ satisfying (3.1). Since there infinitely many such\ $n$\ and\
$\und{n\to\iy}{\text {lim}}\ r_n=1$,\ we come to a contradiction.
\vk.3cm
2)\ \ {\it The case}\ $p>2$:\ \ Since the Jung constant\ $J(\el_p)=\dfrac{1}{\sqrt[p]{2}}$,\ thus\ $d(A)=\sqrt[p]{2}$.\ Writting\ $y_i-c_n=(z^i_1, z^i_2, \ldots),\ \ J(y_i-c_n)=(v^i_1, v^i_2, \ldots)$\ with\ $v^i_j=$sgn$(z^i_j).\|z^i_j\|^{p-1},\ j=1, 2, \ldots$\ and applying part (ii) of Lemma 3.1 coordinatewise to\ $\big\lag (y_i-c_n)-(y_j-c_n),J(y_i-c_n)-J(y_j-c_n)\big\rag$\ we obtain
 
  $$\begin{aligned} \big\lag (y_i-c_n)-(y_j-c_n),J(y_i- & c_n)-J(y_j-c_n)\big\rag 
  =\und{k=1}{\ov{\iy}{\sum}}\ (z^i_k-z^j_k).(v^i_k-v^j_k)\le\\
  &\le \und{k=1}{\ov{\iy}{\sum}}\ \|z^i_k-z^j_k\|^p+
(p-1).\und{k=1}{\ov{\iy}{\sum}}\ (z^i_k.v^j_k+z^j_k.v^i_k),
  \end{aligned}$$

\no and therefore after summing up

  $$\begin{aligned}
r_n^{p-1}.T_n & \le \und{i, j\in I_n}{\sum}\ \al_i.\al_j.\|y_i-y_j\|^p+(p-1).\und{i, j\in I_n}{\sum}\ \al_i.\al_j.\und{k=1}{\ov{\iy}{\sum}}\ (z^i_k.v^j_k+z^j_k.v^i_k)=\\
  & =\und{i, j\in I_n}{\sum}\ \al_i.\al_j.\|y_i-y_j\|^p+2.(p-1).
\und{k=1}{\ov{\iy}{\sum}}\ \Big(\und{i\in I_n}{\sum}\ \al_i.z^i_k\Big).\Big(
\und{j\in I_n}{\sum}\ \al_j.v^j_k\Big)=\\
  & =\und{i, j\in I_n}{\sum}\ \al_i.\al_j.\|y_i-y_j\|^p.
  \end{aligned}$$

The last equality follows from\ $\und{j\in I_n}{\sum}\ \al_j.v^j_k=0,\ k=1, 2, \ldots,$\ which in turn can be deduced from the condition\ $\und{j\in I_n}{\sum}\ \al_j.f_j=0$,\ or equivalently\ $\und{j\in I_n}{\sum}\ \al_j.J(y_j-c_n)=0$.\ Thus in view of (3.2)

  $$\begin{aligned} 2.r_n^p & \le \und{i, j\in I_n}{\sum}\ \al_i.\al_j.\|y_i-y_j\|^p=\\
  & =\und{i, j\in J_n}{\sum}\ \al_i.\al_j.\|y_i-y_j\|^p+
\und{(i,j)\in I^2_n\setminus J^2_n}{\sum}\ \al_i.\al_j.\|y_i-y_j\|^p\le\\
  & \le \big( \sqrt[p]{2}-\varep_0 \big)^p.\Big(\und{i, j\in J_n}{\sum}\ \al_i.\al_j\Big)+2.\Big( 1-\und{i, j\in J_n}{\sum}\ \al_i.\al_j\Big)=
  \end{aligned}$$

 \begin{align}
  & =2-\Big[ 2-\big( \sqrt[p]{2}-\varep_0 \big)^p\Big] .\Big(\und{i, j\in J_n}{\sum}\ \al_i.\al_j\Big)\le\notag\\
  & ({\text {for all $n$ satisfying}}\ (3.1))\notag\\
  & \le 2-\Big[ 2-\big( \sqrt[p]{2}-\varep_0 \big)^p\Big]. \dfrac{1}{m^2},\notag \end{align}
\vk.1cm
\no a contradiction, because\ $\und{n\to\iy}{\text {lim}}\ r^p_n=1$\ and there are infinitely many\ $n$\ satisfying (3.1).
\vk.2cm
One concludes that\ $\al(A_\iy)=d(A)$,\ and hence\ $\al(A)=d(A)$.
\vk.2cm
The proof of Theorem 3.2 is complete.
\vk.5cm
As an immediate consequence one obtains an extension of Gulevich's result for\ $\el_p$\ spaces. 
\vk.5cm
\textbf{Corollary\ 3.3}\ ({\it cf.} \cite{Gul}).\quad {\it Let $A$ be a relatively compact set in an $\el_p$ space with $d(A)>0$. Then

  $$\begin{aligned} & r_{\el_p}(A)<\dfrac{1}{\sqrt[q]{2}}.d(A),\qquad 1<p\le 2;\\
    & r_{\el_p}(A)<\dfrac{1}{\sqrt[p]{2}}.d(A),\qquad 2<p<\iy.
  \end{aligned}$$}
\vk.3cm
\textbf{Theorem\ 3.4.}\quad {\it Let\ $A$\ be an extremal set in a given\ $\el_p$\ space with\ $1<p<\iy$. Then for every\ $\varep\in (0,d(A))$,\ every positive integer\ $m$,\ there exists an\ $m$-simplex\ $\De$\ with vertices in\ $A$\ such that each edge of\ $\De$\ has length not less than\ $d(A)-\varep$.}
\vk.2cm
{\it Proof}.\quad We shall assume\ $r_{\el_p}(A)=1$.\ From the proof of Theorem 3.2 we derived a sequence\ $\{y_{s_{n-1}+1},\ y_{s_{n-1}+2},\ \ldots,\ y_{s_n}\}^\iy_{n=2}$\ in\ $A$,\ a sequence of continuous linear functionals $\{f_{s_{n-1}+1},\ f_{s_{n-1}+2},\ \cdots,\ f_{s_n}\}^\iy_{n=2}$\ in\ $\el_q$ and a sequence of positive numbers\ $\{\al_{s_{n-1}+1},\ \al_{s_{n1}+2},$ \break \ldots,\ $\al_{s_n}\}^\iy_{n=2}$\ (with convention\ $s_1=0$)\ such that:
\vk.2cm
{\rm (i)}\ \ $\lag y_i-c_n,f_i\rag=r_n$\ for\ \ $i\in I_n\col=\{s_{n-1}+1, s_{n-1}+2, \ldots s_n\}$;
\vk.1cm
{\rm (ii)}\ \ $\|f_i\|=1,\ \ i\in I_n$;
\vk.1cm
{\rm (iii)}\ \ $\und{i\in I_n}{\sum}\ \al_i\ f_i=0,\ \ \und{i\in I_n}{\sum}\ \al_i=1$,

\no where\ $c_n\in\el_p$,\ and\ $r_n\in (1-\frac{1}{n},1]$.\ Also we have

 $$\begin{aligned} & 2.r^p_n\le 2^{2-p}.\und{i, j\in I_n}{\sum}\ \al_i.\al_j.\|y_i-y_j\|^p,\qquad 1<p\le 2;\\
    & 2.r^p_n\le \und{i, j\in I_n}{\sum}\ \al_i.\al_j.\|y_i-y_j\|^p,\qquad \qquad 2<p<\iy.
  \end{aligned}$$
\vk.2cm
1)\ \ {\it The case}\ $1<p\le 2$.\quad We denote by

  \begin{align} & T_{nj}\col=2^{2-p}.\und{i\in I_n}{\sum}\ \al_i.\|y_i-y_j\|^p,\notag\\
  & S_n\col=\big\{ j\in I_n:\ \ T_{nj}\ge 2.r^p_n .\big( 1-
 \sqrt{1-r^p_n} \big) \big\},\notag\\
  & \la_n\col=\und{i\in I_n\setminus S_n}{\sum}\ \al_i.\notag
  \end{align}

We have

  $$\begin{aligned} 2.r^p_n & \le 2^{2-p}.\und{i, j\in I_n}{\sum}\ \al_i.\al_j.\|y_i-y_j\|^p=\\
  & =\und{j\in S_n}{\sum}\ \und{i\in I_n}{\sum}\ 2^{2-p}.\al_i.\al_j.\|y_i-y_j\|^p+\und{j\in I_n\setminus S_n}{\sum}\ \und{i\in I_n}{\sum}\ 2^{2-p}.\al_i.\al_j.\|y_i-y_j\|^p\le\\
  & \le 2^{2-p}.2^{{}^{\frac{p}{q}}}.\und{j\in S_n}{\sum}\ \al_j+\und{j\in I_n\setminus S_n}{\sum}\ \al_j.2.r^p_n.\big( 1-\sqrt{1-r^p_n} \big)=\\
  & =2.(1-\la_n)+2.\la_n.r^p_n .\big( 1-\sqrt{1-r^p_n} \big)=\\
  & =2-2.\la_n.\big( 1-r^p_n+r^p_n.\sqrt{1-r^p_n} \big)\le\\
  & \le 2-2.\la_n.\sqrt{1-r^p_n}.
 \end{aligned}$$

Hence

  $$\la_n\le \sqrt{1-r^p_n}\rg 0$$

\no as\ $n\rg\iy$. Thus

  $$\und{n\to\iy}{\text {lim}}\ \Big(\und{i\in S_n}{\sum}\ \al_i\Big)=\und{n\to\iy}{\text {lim}}\ (1-\la_n)=1.$$
\vk.2cm
On the other hand

 $$\begin{aligned} \und{i, j\in I_n}{\sum}\ \al_i.\al_j.\|y_i-y_j\|^p &
   =\und{i\ne j\in I_n}{\sum}\ \al_i.\al_j.\|y_i-y_j\|^p\le\\
   & \le 2^{{}^{\frac{p}{q}}}.\bigg( 1-\Big(\und{i\in I_n}{\sum}\ \al_i^2\Big)\bigg). \end{aligned}$$

So

  $$2.r^p_n\le 2^{2-p}.2^{{}^{\frac{p}{q}}}.\bigg[ 1-\Big(\und{i\in I_n}{\sum}\ \al_i^2\Big)\bigg]\le 2.(1-\al_i^2)$$

\no for each\ $i\in I_n$.\ Therefore\ $\al_i\le \sqrt{1-r^p_n}\rg 0$\ as\
$n\rg\iy$. One concludes that

  $$\big| S_n\big|\lrg\iy\quad {\text {as}}\quad n\rg\iy.$$
\vk.2cm
This implies that for each positive integer\ $m$\ one can choose\ $n$\ sufficiently large so that\ $\big| S_n\big|>m$.\ For $j\in S_n$\ put

  $$\begin{aligned} & S_n(y_j)\col=\bigg\{ i\in I_n:\ \ \|y_i-y_j\|^p\ge
2^{{}^{\frac{p}{q}}}.\Big( 1-\dfrac{1}{\sqrt[4]{n}}\Big)\bigg\},\\
  & \hat{S}_n(y_j)\col=\big\{ y_i:\ \ i\in S_n(y_j) \big\}.
  \end{aligned}$$

As\ $j\in S_n$\ we have

 $$\begin{aligned} 2.\Big( 1-\dfrac{1}{n}\Big)^p. & \bigg( 1-
\sqrt{1-\Big( 1-\dfrac{1}{n}\Big)^p} \bigg) \le 2.r^p_n.\Big( 1-\sqrt{1-r^p_n}\Big)\le\\
 & \le 2^{2-p}.\und{i\in I_n}{\sum}\ \al_i.\|y_i-y_j\|^p=\\
 & = 2^{2-p}.\und{i\in S_n(y_j)}{\sum}\ \al_i.\|y_i-y_j\|^p+
2^{2-p}.\und{i\in I_n\setminus S_n(y_j)}{\sum}\ \al_i.\|y_i-y_j\|^p\le\\
 & \le 2^{2-p}.2^{{}^{\frac{p}{q}}}.\Big( 1-\und{i\in I_n\setminus S_n(y_j)}{\sum}\ \al_i\Big)+2^{2-p}.2^{{}^{\frac{p}{q}}}.\Big( 1-\dfrac{1}{\sqrt[4]{n}}\Big).\Big(\und{i\in I_n\setminus S_n(y_j)}{\sum}\ \al_i\Big)=\\
 & = 2-\dfrac{2}{\sqrt[4]{n}}.\Big(\und{i\in I_n\setminus S_n(y_j)}{\sum}\ \al_i\Big). \end{aligned}$$

Hence

  $$\begin{aligned} \und{i\in I_n\setminus S_n(y_j)}{\sum}\ \al_i & < \bigg[ 1-\Big( 1-\dfrac{1}{n}\Big)^p.\bigg( 1-
\sqrt{1-\Big( 1-\dfrac{1}{n}\Big)^p} \bigg)\bigg].\sqrt[4]{n}=\\
  & =\bigg[ 1-\Big( 1-\dfrac{1}{n}\Big)^p\bigg].\sqrt[4]{n}+ 
 \Big( 1-\dfrac{1}{n}\Big)^p.\sqrt{1-\Big( 1-\dfrac{1}{n}\Big)^p}.\sqrt[4]{n}<\\
  & < \dfrac{p}{\sqrt[4]{n^3}}+\sqrt{\dfrac{p}{\sqrt{n}}}<\\
  & < \dfrac{2.p}{\sqrt[4]{n}}\hk12cm (3.4) 
\end{aligned}$$

\no and so

   $$\und{i\in S_n(y_j)}{\sum}\ \al_i>1-\dfrac{2.p}{\sqrt[4]{n}}\eqno{(3.5)}$$
\vk.2cm
Now for a given positive integer\ $m$\ we choose\ $n$\ sufficiently large such that

  $$\big| S_n\big|>m\qquad {\text {and}}\qquad \dfrac{2.p.m}{\sqrt[4]{n}}<1.$$
\vk.2cm
We claim that for every choice of\ $i_1,\ i_2,\ \ldots,\ i_m\in S_n$\ we have

  $$\und{k=1}{\ov{m}{\bigcap}}\ S_n(y_{i_k})\ne\emptyset\eqno{(3.6)}$$

Indeed, otherwise\ $\und{k=1}{\ov{m}{\bigcap}}\ S_n(y_{i_k})=\emptyset$\ would imply 

  $$S_n(y_{i_1})\subset I_n\setminus \Big( \und{k=2}{\ov{m}{\bigcap}}\ S_n(y_{i_k})\Big)=\und{k=2}{\ov{m}{\bigcup}}\ \Big( I_n\setminus S_n(y_{i_k})\Big).$$

\no Consequently by (3.4) and (3.5)

  $$\begin{aligned} 1-\dfrac{2.p}{\sqrt[4]{n}} & < \und{\nu\in S_n(y_{i_1})}{\sum}\ \al_\nu\le\\
  & \le \und{k=2}{\ov{m}{\sum}}\ \und{\nu\in I_n\setminus S_n(y_{i_k})}{\sum}\ \al_\nu < (m-1).\dfrac{2.p}{\sqrt[4]{n}},\end{aligned}$$

\no or\ $1<\dfrac{2.p.m}{\sqrt[4]{n}}$,\ a contradiction to the choice of\ $n$.
\vk.3cm
Furthermore from (3.6) it follows that if\ $1\le k\le m$\ and\ $i_1,\ i_2,\ \ldots,\ i_k\in S_n$,\ then
  
   $$\und{\nu=1}{\ov{k}{\bigcap}}\ \hat{S}_n(y_{i_\nu})\ne\emptyset.$$
\vk.2cm
With\ $m$\ and\ $n$\ as above let's fix some\ $j\in S_n$.\ Setting\ $z_1\col=y_j$\ we take consecutively\ $z_2\in \hat{S}_n(z_1),\ z_3\in \hat{S}_n(z_1)\cap \hat{S}_n(z_2),\ \ldots,\ z_{m+1}\in \und{k=1}{\ov{m}{\bigcap}}\ \hat{S}_n(z_k)$. Obviously

  $$\|z_i-z_j\|^p\ge 2^{{}^{\frac{p}{q}}}.\Big( 1-\dfrac{1}{\sqrt[4]{n}}\Big)$$

\no for all\ $i\ne j$\ in\ $\{1, 2, \ldots, m+1\}$.\ Now for a given\ $\varep\in (0,\sqrt[q]{2})$\ one can choose\ $n$\ sufficiently large as above, and moreover so that

  $$2^{{}^{\frac{p}{q}}}.\Big( 1-\dfrac{1}{\sqrt[4]{n}}\Big)\ge \Big( \sqrt[q]{2}-\varep\Big)^p.$$
\vk.2cm
One obtains an\ $m$-simplex\ formed by\ $z_1,\ z_2,\ \ldots,\ z_{m+1}$,\ whose edges have length not less than\ $\Big( \sqrt[q]{2}-\varep\Big)$,\ as claimed.
\vk.4cm
2)\ \ {\it The case}\ $2<p<\iy$\ can be proceeded in the same way just with replacing those\ $T_{nj}$\ and\ $S_n$\ in the first case suitably, {\it i.e.} as 

  \begin{align} & T_{nj}\col=\und{i\in I_n}{\sum}\ \al_i.\|y_i-y_j\|^p,\notag\\
  & S_n\col=\big\{ j\in I_n:\ \ T_{nj}\ge 2.r^p_n .\big( 1-
 \sqrt{1-r^p_n} \big) \big\}.\notag\end{align}
\vk.2cm
The proof of Theorem 3.4 is complete.
\vk.5cm

\end{document}